# A multi-objective sustainable planning for a real hazardous waste production problem

Abed Zabihian-Bisheh[1], Hadi Rezaei Vandchali[2], Vahid Kayvanfar[3], Frank Werner[4,*]

[1] *School of Industrial Engineering, Collage of Engineering, University of Tehran, Tehran, Iran*

[2] *Australian Maritime College, University of Tasmania, Launceston, Australia*

[3] *Department of Industrial Engineering, Sharif University of Technology, Tehran, Iran*

[4] *Faculty of Mathematics, Otto-von-Guericke-University, Magdeburg, Germany*

**Abstract**

A significant amount of hazardous waste generated from health sectors and industrial processes has posed a major threat to human health by causing environmental issues and contamination of air, soil, and water resources. This paper presents a multi-objective mixed-integer nonlinear programming (MINLP) formulation for a sustainable hazardous waste location-routing problem. The location of the facilities and routing decisions for transporting hazardous waste and the waste residue is considered to design a suitable waste collection system. The presented model simultaneously minimizes the total costs of the waste management system, total risks from transportation and facilities, along with $CO_2$ emissions. A real-world case study is presented to illustrate the applicability of the proposed model. To illustrate the significance of sustainability, the results of the original model are compared with the results of the model without considering sustainability. It indicates that, under the condition when sustainability is not taken into account, total cost, transportation, and site risk along with $CO_2$ emission increased, which in turn demonstrated the importance of sustainability. Furthermore, the managerial insights gained from the optimal results would enable the managers to make better decisions in the hazardous waste management system.

**Keywords:** Sustainability; Waste management; Location-routing problem; Capacity planning; Case study.

## 1. Introduction

Hazardous waste that is produced by the industrial process can cause several impacts on humans, animals, and plants. Toxicity, reactivity, ignitability, and corrosiveness are four kinds of waste categories, and if any waste possesses at least one of these categories, it would be categorized as hazardous waste (Alumur and Kara 2007). In recent years, population growth as well as industrial and technological developments have increased environmental impacts on our planet, which have raised governmental concerns regarding a hazardous waste management (HWM) system. An HWM system involves a multitude of tasks such as collecting, transporting, recycling, treatment, and disposal of hazardous waste as well as determining related routes for collecting hazardous waste and waste residue to deliver them to facilities (Samanlioglu 2013). Due to various categories of physical and chemically hazardous waste, selecting suitable locations for treatment, recycling, and disposal facilities as well as how to allocate this hazardous

---

* Corresponding Author:

E-mail addresses: abed.zabihian1992@gmail.com (A. Zabihian-Bisheh), hadi.rezaei@utas.edu.au (H.R. Vandchali); kayvanfar@sharif.edu (V. Kayvanfar), frank.werner@ovgu.de (Frank Werner)

waste to the associated facilities plays a crucial role in the HWM system and could be a challenging task (Nema and Gupta 1999). To design an efficient HWM system, several perspectives, including environmental, economic, and social aspects, should be considered simultaneously, making it a more complex problem (Jacobs and Warmerdam 1994).

A sustainable hazardous waste management (SHWM) system is of high importance to guarantee human health and safety. In a waste management system (WMS), sustainability must be addressed properly; otherwise, it fails (Asase et al. 2009), having examples of composting facilities in Delhi (Talyan et al. 2008) and composting and incinerator facilities in Turkey (Kanat 2010). These facilities were opened to avoid the adverse effects of solid waste on the environment that is placed in the open dumps area. However, they failed because they were not efficient and also had high operational costs. The lack of considering sustainability issues during the design phase was the main cause of the failure (Kanat 2010). Incidents like these two failures indicate the significance of sustainability in a waste management system. In this regard, an SHWM must be effective in the environmental perspective, affordable in the economic aspect, and acceptable in social terms (McDougall et al. 2008). The World Summit of Sustainable Development (WSSD) defines sustainability as a balance between environmental protection, economic profit, and social developments, which highlights the importance of environmental protection as one of the pillars of sustainability (White and Lee 2009). To address this issue, this paper considers waste residues and transportation systems as the two sources of emissions of greenhouse gases and aims to minimize the $CO_2$ emissions from waste, waste residues, and transportation systems.

Social responsibility is another important pillar of sustainability on which one should be focused. Due to growing social responsibility concerns, many scientists have tried to deal with environmental risks (De Buck et al. 1999), (Qiu et al. 2001). Numerous studies have been conducted for waste management systems considering risks such as population risk and transportation risk in their models as an objective required to be minimized (Current and Ratick 1995), (Samanlioglu 2013), (Zografros and Samara 1989) (Yu et al. 2020a). However, few studies have considered the impact of social responsibility on their proposed models (Garrido and Bronfman 2017). This paper addresses the social responsibility issue by considering the transportation risk and site risk as an objective function to be minimized.

To prevent interaction between different types of incompatible waste, they should be carried separately, including emissions and the generation of heat. In this regard, an inhomogeneous fleet of vehicles is exploited to collect the waste in a divided fleet that is compatible with the load. The network in our problem consists of five components, including generation nodes, a central depot, and three facilities (recycling, treatment, and disposal). In addition, there is a capacity level for each facility that allows us to utilize these facilities effectively. The capacity level includes three levels for each facility, which are determined by the amount of waste delivered to the facilities. In this paper, we propose an MINLP model, in which the location of the facilities and routing decisions are considered simultaneously in a sustainable hazardous waste location-routing problems (SHWLRPs). In this regard, a new mathematical programming model is developed to minimize the total investment cost, transporting cost and $CO_2$ emissions from the transportation system and waste residues. In addition, we aim to minimize transportation and site risks.

The remainder of this study is as follows: Section 2 presents a review of the related literature. The descriptions and mathematical formulations are presented in Section 3. In Section 4, the importance and the value of the proposed model are evaluated through a real case study. Finally, Section 5 provides some conclusions and suggests some future directions to further advance the paper.

## 2. Literature review

Three main objective functions are generally addressed by the researchers to consider simultaneously different aspects of the undesirable facility location problem in the hazardous waste management; the minimization of costs including construction and operational costs, the minimization of transportation and facilities risks which the population around those geographic areas are being exposed to, and the maximization of the equity of the distribution risk which is obtained by the maximum zonal risk per unit. These three objectives were first developed by (Ratick and White 1988) in the facility location problem. Later, many researchers assessed the tradeoffs between risk, cost, and the distribution of equitable risk simultaneously (Erkut and Verter 1998), (Wyman and Kuby 1995), (Zhang and Zhao 2011), (Asgari et al. 2017), (Zhao and Huang 2019). Since hazardous waste is extremely harmful to our environment and human health, it attracts a public and governmental concern, which encourages many researchers to develop a mathematical model approach to help the decision-makers in evaluating hazardous waste. Two significant aspects of dealing with hazardous waste are the location of facilities and the routing of hazardous waste between the facilities (Current and Ratick 1995). Regarding the first aspect, the primary work associated with a semi-desirable or partially noxious facility location problem was developed by (Goldman and Dearing 1975) in which they proposed the notion of an obnoxious facility on a network. (Erkut and Neuma 1989) surveyed location models in the literature of operations research. Later, concerns about the location of undesirable facilities have greatly raised as the magnitude of hazardous materials increased. For example, (Melachrionoudis et al. 1995) developed a multi-objective mixed-integer programming model for sitting landfills. Minimizing the total cost, the risk to nearby population centers, risk nearby to the ecosystem, and inequity of risk distribution were four objectives that they considered in their paper. (Tuzkaya et al. 2008) addressed the undesirable facility location selection problem to determine an appropriate location to construct undesirable facilities in Istanbul. Benefits, opportunities, costs, and risks were taken into account as criteria to find optimal locations. (Darmian et al. 2020) presented a sustainable multi-objective location-districting optimization model for collecting municipal solid waste. The purpose in this paper was the design of an efficient system for municipal services by combing the decisions with regard to urban area districts and the location of waste collection centers.

Regarding the second aspect, (Minh et al. 2013) proposed a multi-objective vehicle routing problem for collecting hazardous waste. They applied a memetic algorithm to minimize the total traveling time and the number of vehicles for collecting the waste. (Das and Bhattacharyya 2015) proposed a MILP model to minimize the length of the routes for collecting the waste to reduce the total costs of the system in India. (Louati 2016), to collect municipal solid waste, proposed a generalized vehicle routing model, including heterogeneous vehicles within time windows, and they applied a MILP approach to reduce both total traveling distances and operational hours of the vehicles. (Bronfman et al. 2016) developed a new method to the problem of hazardous material transportation in urban areas in addition to the minimization of general population risk. (Zahiri et al. 2020) presented a transportation network for hazardous material. To this aim, they developed a new bi-objective mathematical model to design the hazardous material transportation network.

Considering the location and routing decisions in HWM leads to a Hazardous Waste Location-Routing Problem (HWLRP). The first HWLRP was proposed by (Zografros and Samara 1989). They considered one type of hazardous waste and three objectives, including travel time, transportation risk, and site risk. They proposed a goal programming model for determining the location and routing decisions. (ReVelle et al. 1991) proposed a location-routing problem to minimize the joint facility location, transportation route, cost, and perceived risks. (Stowers and Palekar 1993) proposed a location model for an obnoxious facility location problem in which they simultaneously incorporated routing decisions into their model. (Ghaderi and Burdett 2019) considered a bi-modal transportation network including road and rail links. They formulated a two-stage stocashtic programming model to minimize the transportation cost and risk. (Rabbani et al. 2020) proposed a mixed-integer bi-objective waste management location-routing problem in an automotive company in Iran. For this purpose, they considered a heterogeneous fleet of vehicles to determine an optimal route for collecting the waste and the number of vehicles in the transportation phase. (Rabbani et al. 2019) incorporated a multi-period HWLRP. They applied a multi-objective MINLP model to select the best route for collecting the waste. (Yu et al. 2020b) developed a two-echelon multi-objective HWLRP. Finding the best route of the vehicles for collecting hazardous waste to minimize cost and risk was one of the leading purposes of their study. (Mohri et al. 2020) addressed a hazmat routing-scheduling problem to minimize risk in hazmat land transportation by simultaneously incorporating routing and scheduling approaches. The routing approach simultaneously divides the routes of the transported hazmat while the scheduling approach divides the transportation time slots on the same routes.

Due to various types of hazardous waste in the real world, many works have tried to address those types differently. For example, (Samanlioglu 2013) proposed different kinds of an industrial hazardous waste network design. This paper presented a new multi-objective model to minimize total cost, transportation, and site simultaneously. (Nema and Gupta 2003) and (Nema and Gupta 1999) are two other works that addressed a location-routing model to consider different categories of hazardous waste. They suggested waste-waste and waste-technology compatibility constraints and presented a multi-objective goal programming model to find the location of the disposal and treatment centers as well as transportation routes to these facilities. (Ghezavati and Morakabchian 2015) and (Rabbani et al. 2018) described an HWLRP with recycling, treatment, and disposal facilities. (Rabbani et al. 2018) addressed waste-waste compatibility to guarantee that a type of waste is not collected together with incompatible waste, and all of them are collected with an inhomogeneous fleet of vehicles. Furthermore, (Ghezavati and Morakabchian 2015) divided the generation nodes into two distinctive types of centers. The first type is where the waste is produced, and the second type is where the waste is collected. (Aydemir-Karadag 2018) formulated a mathematical model for a profit-HWLRP in which multiple types of waste were taken into account. They considered two treatment technologies, including chemical treatment and incineration to deal with different categories of waste delivered to the facilities in accordance with their classification.

As discussed earlier, the environmental aspect is one of the essential pillars of sustainability. The environmental perspective was rarely taken into account in the HWM literature. (Chang et al. 2012) formulated a solid waste management problem in which the global warming and cost-profit factor were considered to obtain an optimal planning for the waste system. They considered a scenario-based design process to estimate greenhouse gas emissions. (Minoglou and Komilis 2013) proposed an integrated solid waste management with regard to two objective functions, including total cost minimization and the minimization of the carbon dioxide ($CO_2$) emissions.

Table 1 presents a comparative analysis of the published studies that were particularly focusing on the field of HWLRPs from 1989 to 2020. According to this table and the reviewed studies in this section, it has been found that little attempts have been made to consider all the three components (economic, environmental and social) simultaneously in the HWLRP literature. In addition, to the best of our knowledge, capacity planning of the facilities has not been taken into consideration in the HWLRP area. However, due to continuous changes in the relationship between waste generation and the facility capacities' availability, it is important to have a capacity planning for the effective utilization of the facilities (Huang et al. 1997). To overcome these shortcomings and fulfill these gaps, this paper develops a model which considers three aspects of sustainability at the same time and also formulates the capacity planning for each facility. The main contributions of this paper, which distinguish our work from the previous studies in this area, are as follows:
- Developing a sustainable location-routing hazardous waste management system with capacity planning for each facility;
- Proposing a new environmental objective function to minimize the environmental impact of $CO_2$ emissions in the location-routing hazardous waste management;
- Considering a heterogeneous fleet of vehicles to prevent interactions between different waste types.

## 3. Problem description

In this section, first, the description of the proposed HWLRP and the assumptions characterized by the given problem are presented. Then, the mathematical model of the resulting HWLRP is formulated using an MINLP approach under the foregoing assumptions. Finally, the augmented ε-constraint method is applied to convert the HWLRP to a single-objective problem.

*3.1. Description*

Fig. 1 presents the schematic view of the structure of the concerned HWLRP and the links between the facilities and generation nodes. The network of the developed HWRLP includes the generation nodes, recycling, disposal, and treatment facilities, and a central depot which is the starting and ending point of each route. Moreover, there are potential locations for each set of facilities and a set of existing facilities in the network. Hazardous waste is produced at a generation node, and it can be hospitals, factories, and collection centers. There are many potential nodes, where the recycling facilities and treatment facilities can be established. A link that connects the two adjacent facilities shows the physical roads in which the waste is transported between the facilities and generation nodes through these roads.

The presented model can handle different kinds of waste that are incompatible with each other. Four types of waste have been considered in this paper. These four types are (I) recyclable waste, (II) non-recyclable waste, but it is suitable for incineration technology, (III) non-recyclable waste, but it is suitable for chemical technology, and (IV) waste that is non-recyclable, but it is compatible with both technologies. An inhomogeneous fleet of vehicles is considered to collect the incompatible types of waste in different vehicles to avoid interactions between incompatible types of waste. There is a limitation for the vehicle that carries a specific type of waste associated with its capacity and the length of the tour it travels. In this location-routing problem for collecting the waste, each vehicle starts from the central depot, and after unloading the waste, it comes back to the central depot. Hazardous waste is accumulated at the generation nodes. The accumulated waste is collected with a compatible vehicle, and then the recyclable amounts are shipped to the recycling facilities, and the non-recyclable amounts are transferred to the treatment facility. Again, at the treatment facility, the amount of waste residue that is recyclable is delivered to recycling facilities, and the remaining ones which are non-recyclable are sent to disposal facilities. Also, waste residues at recycling facilities are routed to the disposal facility.

We optimize the SHWLRP by considering the minimization of total costs, total risks, and $CO_2$ emissions. The total costs include the total transportation cost and the total investment cost of opening facilities. The second objective deals with reducing the total risks. The total risks are composed of the transportation risk and the site risk. Exposure of the population within the transportation routes is defined as transportation risk. The site risk is defined as the risk that the population around the treatment and disposal facilities are facing. This risk is measured by the amount of loads at each facility and the population living around the facilities within a given radius.

While addressing the economic and social aspects of the system, to design a sustainable hazardous waste management system, the environmental aspect should be incorporated as well. In this regard, we design an environmental-friendly HWLRP by minimizing the $CO_2$ emissions.

The main assumptions to facilitate the model formulation are summarized as follows:

- Parameters are considered to be deterministic.
- All the recycling, treatment, and disposal facilities have limited capacities.
- Each generation node is served only one time by one vehicle for each type of waste.
- The number of potential facilities that can be opened is limited.
- The transportation cost for the vehicles is related to the distance they traveled.
- Two kinds of treatment technology exist, and at each treatment center, at most one of them is allowed to be established.

**Table 1**
Most relevant literature on the hazardous waste location-routing problem.

| Study | Decision variables | | | Objective functions | | | Vehicles | | Compatibility | |
|---|---|---|---|---|---|---|---|---|---|---|
| | Location | Allocation/ Routing | Capacity | Total cost | Risk | $CO_2$ emission | Homogenous | Heterogeneous | Waste-Waste | Waste-Technology |
| Zografros and Samara (1989) | ✓ | ✓ | | ✓ | ✓ | | ✓ | | | |
| List and Mirchandani (1991) | ✓ | ✓ | | ✓ | ✓ | | | | | |
| Stowers and Palekar (1993) | ✓ | ✓ | | | ✓ | | ✓ | | | |
| Current and Ratick (1995) | ✓ | ✓ | | ✓ | ✓ | | ✓ | | | |
| Nema and Gupta (2003) | ✓ | ✓ | | ✓ | ✓ | | ✓ | | ✓ | ✓ |
| Alumur and Kara (2007) | ✓ | ✓ | | ✓ | ✓ | | ✓ | | | ✓ |
| Zhang and Zhao (2011) | ✓ | ✓ | | ✓ | ✓ | | ✓ | | | |
| Samanlioglu (2013) | ✓ | ✓ | | ✓ | ✓ | | ✓ | | | ✓ |
| Jiahong Zhao and Verter (2015) | ✓ | ✓ | | ✓ | ✓ | | ✓ | | | ✓ |
| Ghezavati and Morakabchian (2015) | ✓ | ✓ | | ✓ | ✓ | | ✓ | | | |
| Zhao et al. (2016) | ✓ | ✓ | | ✓ | ✓ | | ✓ | | | |
| Asgari et al. (2017) | ✓ | ✓ | | ✓ | ✓ | | ✓ | | | |
| Yilmaz et al. (2017) | ✓ | ✓ | | | | | ✓ | | | ✓ |
| Farrokhi-Asl et al. (2017) | ✓ | ✓ | | | | | | ✓ | | |
| Rabbani et al. (2018) | ✓ | ✓ | | ✓ | ✓ | | | ✓ | ✓ | ✓ |
| Aydemir-Karadag (2018) | ✓ | ✓ | | ✓ | | | | | | ✓ |
| Zhao and Huang (2019) | ✓ | ✓ | | | | | ✓ | | | |
| Rabbani et al. (2019) | ✓ | ✓ | | ✓ | ✓ | | | ✓ | | |
| Yu et al. (2020a) | ✓ | ✓ | | ✓ | ✓ | | ✓ | | | |
| Yu et al. (2020b) | ✓ | ✓ | | ✓ | ✓ | | | | | |
| This study | ✓ | ✓ | ✓ | ✓ | ✓ | ✓ | | ✓ | ✓ | ✓ |

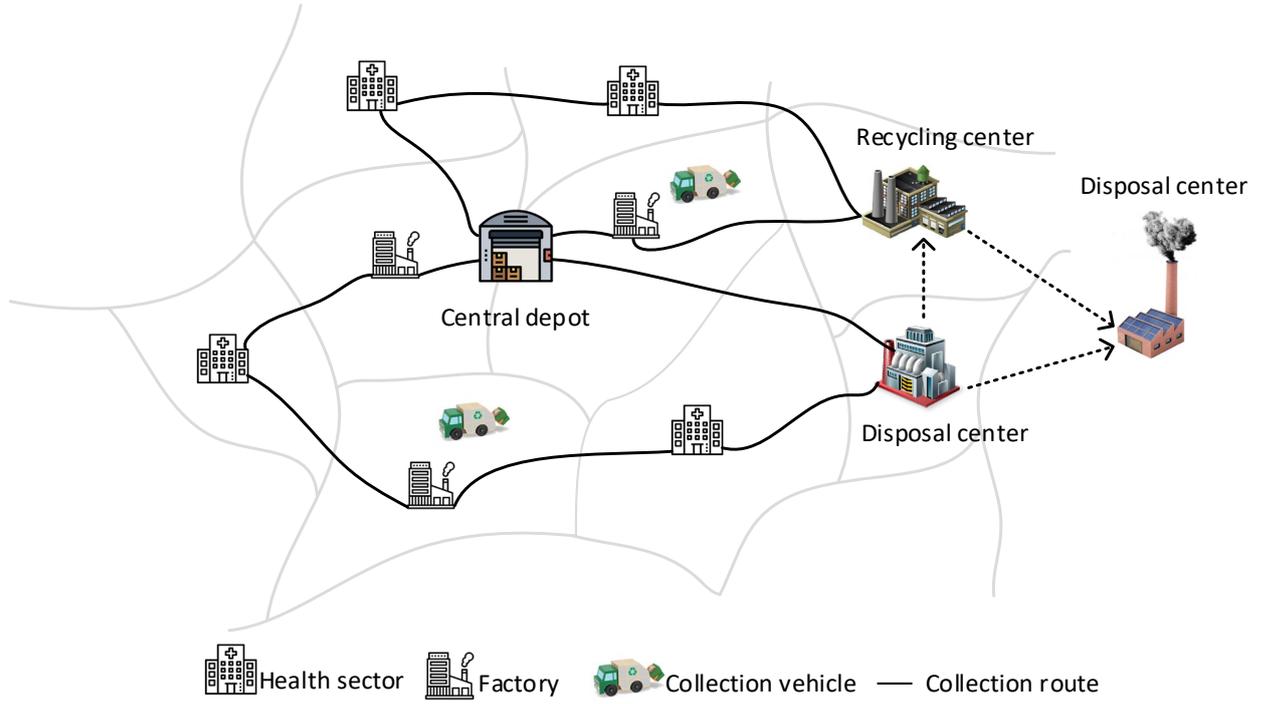

**Fig. 1.** The framework of the hazardous waste management system.

*3.2. Model formulation*

We propose a multi-objective MINLP for HWLRs, which simultaneously determines the locations of the facilities, the routing decisions, the amounts of loads that are operated at each facility. The sets, indices, parameters, and decision variables of the model are as follows:

Sets:
| | |
|---|---|
| $N$ | Set of nodes in the network |
| $G$ | Set of generation nodes $G = 1,2,\ldots,g$ |
| $R$ | Set of recycling centers $R = 1,2,\ldots,r$ |
| $\acute{R}$ | Set of existing recycling centers $\acute{R} \in R$ |
| $T$ | Set of treatment centers $T = 1,2,\ldots,t$ |
| $\acute{T}$ | Set of existing treatment centers $\acute{T} \in T$ |
| $D$ | Set of disposal centers $D = 1,2,\ldots,d$ |
| $\acute{D}$ | Set of existing disposal centers $\acute{D} \in D$ |
| $W$ | Set of hazardous waste types $W = 1,2,\ldots,w$ |
| $Q$ | Set of treatment technologies $Q = 1,2,\ldots,q$ |
| $F$ | Set of depots $F = 1,2,\ldots,f$ |
| $K$ | Set of the collection fleet $K = 1,2,\ldots,k$ |
| $H$ | Set of available capacity levels for establishing the facilities $H = 1,2,\ldots,h$ |

Parameters:
| | |
|---|---|
| $c_{ij}$ | Transportation cost for one unit of waste on link $(i,j) \in N$, |
| $ft_{qih}$ | Investment cost of establishing a treatment facility with technology $q$ and a capacity level of $h$ at node $i \in T$ |
| $fr_{ih}$ | Investment cost of establishing a recycling facility with a capacity level of $h$ at node $i \in R$ |
| $fd_{ih}$ | Investment cost of establishing a disposal facility with a capacity level of $h$ at node $i \in D$ |
| $d_{wi}$ | Quantity of waste type $w \in W$ accumulated at generation node $i \in G$ |
| $Ptr_{ij}$ | Transportation risk per one kg of hazardous waste from node $i$ to node $j$, $i \in T$, $j \in R$ |
| $Ptd_{ij}$ | Transportation risk per one kg of hazardous waste from node $i$ to node $j$, $i \in T$, $j \in D$ |
| $Prd_{ij}$ | Transportation risk per one kg of hazardous waste from node $i$ to node $j$, $i \in R$, $j \in D$ |
| $s_{rh}$ | Operating risk of a recycling center with a capacity level of $h$ at candidate node $r \in R$ |
| $s_{thq}$ | Operating risk of a treatment center with a capacity level of $h$ at node $t \in T$ with technology $q \in Q$ |
| $s_{dh}$ | Operating risk of a disposal center with capacity level $h$ at candidate node $d \in D$ |
| $Pc_{ij}$ | Maximum allowable risk tolerance capacity on link $(i.j) \in N$ |
| $\beta_{wq}$ | Ratio of recyclable hazardous waste type $w$ which are treated with technology $q \in Q$ |
| $r_{wq}$ | Ratio of reduction of mass for hazardous waste type $w$ which are treated with technology $q \in Q$ |
| $\gamma_i$ | Ratio of total waste that is recycled at node $i \in R$ |
| $tc_{qjh}$ | Capacity of treatment technology $q \in Q$ with capacity level $h$ at node $j \in T$ |
| $tc_{qj}^m$ | Minimum quantity of waste required for establishing a treatment technology $q \in Q$ at node $j \in T$ |
| $rc_{jh}$ | Maximum capacity of a recycling center with capacity level $h$ at node $j \in R$ |
| $rc_j^m$ | Minimum quantity of waste required for establishing a recycling facility at node $j \in R$ |
| $dc_{jh}$ | Maximum capacity of a recycling center with capacity level $h$ at node $j \in D$ |
| $dc_j^m$ | Minimum amount of waste residues required for opening a disposal facility at node $j \in D$ |
| $com_{wq}$ | Waste compatibility with technology $q \in Q$; 1 if compatible; 0 otherwise |
| $Qr_{wi}$ | $CO_2$ emissions per one $kg$ of waste type $w \in W$ to be recycled in recycling center $i \in R$ |
| $Qt_{wqi}$ | $CO_2$ emissions per one $kg$ of waste type $w \in W$ to be treated in the treatment center with technology $q \in Q$, $i \in T$ |
| $Qd_{wi}$ | $CO_2$ emissions per one $kg$ of disposed waste type $w \in W$ in node $i \in D$ |
| $QTr_{wij}$ | $CO_2$ emissions per one $km$ and one $kg$ from the transportation of waste type $w \in W$ from node $i$ to node $j$, $i \in T$, $j \in R$ |
| $QTd_{wij}$ | $CO_2$ emissions per one $km$ and one $kg$ from the transportation of waste type $w \in W$ from node $i$ to node $j$, $i \in T, R$, $j \in D$ |
| $dis_{ijk}$ | Travel distance of vehicle $k \in K$ for delivering waste on link $(i,j) \in N$ |
| $ve_{wk}$ | Waste compatibility with vehicle $k \in K$; 1 if compatible; 0 otherwise |
| $a_{qi}$ | Availability of treatment technology $q \in Q$ at the existing treatment facility $i \in \hat{T}$; 1 if available; 0 otherwise |
| $\delta_w$ | Capacity of a vehicle compatible with waste type $w \in W$ |
| $\mu_w$ | Maximum traveling distance of a vehicle compatible with waste type $w \in W$ |

Variables:

| | |
|---|---|
| $x_{ijk}$ | 1 if vehicle $k \in K$ visited node j just after node i; 0 otherwise |
| $z_{ij}$ | Quantity of waste residues delivered from node $i$ to node $j$ $i \in T$, $j \in D$ |
| $k_{ij}$ | Quantity of recyclable waste residue delivered from node $i$ to node $j$ $i \in T$, $j \in R$ |
| $v_{ij}$ | Quantity of waste residues transported from node $i$ to node $j$ $i \in R$, $j \in D$ |
| $xr_i$ | Quantity of waste recycled at node $i \in R$ |
| $xt_{wi}$ | Quantity of hazardous waste type $w \in W$ treated at node $i \in T$ |
| $xd_i$ | Quantity of waste residue disposed at node $i \in D$ |
| $e_{ik}$ | Traveled distance by vehicle $k \in K$ after node $i$ |
| $lo_{ik}$ | Vehicle's load $k \in K$ after node $i$ |
| $r_{ih}$ | 1 if a recycling facility with capacity level $h$ is located at candidate node $i \in R$; 0 otherwise |
| $t_{qih}$ | 1 if a treatment facility with capacity level $h$ is located at candidate node $i \in T$ with technology $q \in Q$; 0 otherwise |
| $d_{ih}$ | 1 if a disposal facility with capacity level $h$ is located at candidate node $i \in D$; 0 otherwise |

### 3.2.2. Proposed mathematical model

The formulation of a mathematical programming model for HWLRPs with respect to the aforementioned notations, which extends the mathematical model proposed by (Rabbani et al. 2018), is as follows:

$$Min f_1(x) = \sum_{i \in G} \sum_{j \in T \cup R \cup D} \sum_{k \in K} c_{ij} x_{ijk} lo_{ik} + \sum_{i \in T} \sum_{j \in D} c_{ij} z_{ij} + \sum_{i \in R} \sum_{j \in D} c_{ij} v_{ij} + \sum_{i \in T} \sum_{j \in R} c_{ij} k_{ij}$$
$$+ \sum_{q \in Q} \sum_{i \in T-T'} \sum_{h \in H} ft_{qih} t_{qih} + \sum_{i \in D-D'} \sum_{h \in H} fd_{ih} d_{ih} + \sum_{i \in R-R'} \sum_{h \in H} fr_{ih} r_{ih} \quad (1)$$

$$Min f_2(x) = \sum_{i \in T} \sum_{j \in R} Ptr_{ij} k_{ij} + \sum_{i \in T} \sum_{j \in D} Ptd_{ij} z_{ij} + \sum_{i \in R} \sum_{j \in D} Prd_{ij} v_{ij} + \sum_{w \in W} \sum_{(i,t) \in T} \sum_{h \in H} s_{th} xt_{wi}$$
$$+ \sum_{(r,i) \in R} \sum_{h \in H} s_{rh} xr_i + \sum_{(d,i) \in D} \sum_{h \in H} s_{dh} xd_i \quad (2)$$

$$Min f_3(x) = \sum_{w \in W} \sum_{i \in R} Qr_{wi} xr_i + \sum_{w \in W} \sum_{i \in R} \sum_{q \in Q} Qt_{wqi} xt_{wi} + \sum_{w \in W} \sum_{i \in D} Qd_{wi} xd_i$$
$$+ \sum_{w \in W} \sum_{i \in T} \sum_{j \in R} QTr_{wij} dis_{ij} k_{ij} + \sum_{w \in W} \sum_{i \in T, R} \sum_{j \in D} QTd_{wij} dis_{ij} (z_{ij} + v_{ij}) \quad (3)$$

Subject to:

$$\sum_{i \in F} \sum_{j \in G} x_{ijk} = 1 \quad \forall k \in K \quad (4)$$

$$\sum_{i \in F \cup G} x_{ijk} - \sum_{i' \in G \cup R \cup T} x_{ji'k} = 0 \quad \forall j \in G, k \in K \quad (5)$$

$$\sum_{j \in G \cup R \cup T} \sum_{k \in K} x_{ijk} ve_{wk} = 1 \quad \forall i \in G, w \in W \quad (6)$$

$$\sum_{i \in G} x_{ijk} - \sum_{i' \in F} x_{ji'k} = 0 \quad \forall j \in R \cup T, k \in K \quad (7)$$

$$x_{ijk} \le \sum_{w \in W} \sum_{q \in Q} \sum_{h \in H} ve_{wk} com_{wk} t_{qih} \quad \forall i \in G, j \in T, k \in K \quad (8)$$

$$x_{ijk} \le \sum_{w \in W} ve_{wk} \frac{\left(\left(2 - \sum_{q \in Q} com_{wq}\right) r_j\right)}{2} \quad \forall i \in G, j \in R, k \in K \quad (9)$$

$$e_{ik} - e_{jk} + \sum_{w \in W} ve_{wk} \left((\mu_w + dis_{ij}) x_{ijk} + (\mu_w - dis_{ij}) x_{jik}\right) \le \sum_{w \in W} ve_{wk} \mu_w \quad \forall i, j \in G \cup R \cup T \quad (10)$$

$$\sum_{i \in F} dis_{ij} x_{ijk} \le e_{jk} \le \sum_{w \in W} ve_{wk} \left(\mu_w + \sum_{i \in F} (dis_{ij} - \mu_w) x_{ijk}\right) \quad \forall j \in G, k \in K \quad (11)$$

$$e_{ik} \le \sum_{w \in W} ve_{wk} \mu_w - \sum_{j \in F} dis_{ij} x_{ijk} \quad \forall j \in R \cup T, k \in K \quad (12)$$

$$lo_{ik} - lo_{jk} + \sum_{w \in W} ve_{wk} \delta_w x_{ijk} \le \sum_{w \in W} ve_{wk} (\delta_w - d_{jw}) \quad \forall i, j \in G, k \in K \quad (13)$$

$$\sum_{w \in W} d_{wi} ve_{wk} \le lo_{ik} \le \sum_{w \in W} ve_{wk} \delta_w \quad \forall i \in G, k \in K \quad (14)$$

$$\sum_{i \in F} \sum_{w \in W} x_{ijk} d_{wj} ve_{wk} \le lo_{jk} \quad \forall j \in G, k \in K \quad (15)$$

$$lo_{jk} \leq \sum_{w \in W} ve_{wk}\left(\delta_w + \sum_{i \in F}(d_{wj} - \delta_w)x_{ijk}\right) \quad \forall j \in G, k \in K \tag{16}$$

$$xt_{wj} = \sum_{k \in K}\sum_{i \in G} x_{ijk} lo_{ik} ve_{wk} \quad \forall j \in T, w \in W \tag{17}$$

$$\sum_{w \in W} xt_{wj} \leq \sum_{q \in Q}\sum_{h \in H} tc_{jh} t_{qjh} \quad \forall j \in T \tag{18}$$

$$\sum_{w \in W} xt_{wj} \geq \sum_{q \in Q}\sum_{h \in H} tc_{jh}^m t_{qjh} \quad \forall j \in T \tag{19}$$

$$xr_j = \sum_{w \in W}\sum_{k \in K}\sum_{i \in G} x_{ijk} lo_{ik} ve_{wk} + \sum_{i' \in T} k_{i'j} \quad \forall j \in R \tag{20}$$

$$\sum_{w \in W}\sum_{q \in Q} xt_{wi}(1 - r_{wq})\beta_{wq} = \sum_{j \in R} k_{ij} \quad \forall i \in T \tag{21}$$

$$xr_j \leq \sum_{h \in H} rc_{jh} r_{jh} \quad \forall j \in R \tag{22}$$

$$xr_j \geq rc_j^m \sum_{h \in H} r_{jh} \quad \forall j \in R \tag{23}$$

$$\sum_{w \in W}\sum_{q \in Q}\sum_{h \in H} xt_{wi} t_{qih}(1 - r_{wq})(1 - \beta_{wq}) = \sum_{j \in D} z_{ij} \quad \forall i \in T \tag{24}$$

$$xr_i(1 - \gamma_i) = \sum_{j \in D} v_{ij} \quad \forall i \in R \tag{25}$$

$$xd_i = \sum_{j \in T} z_{ji} + \sum_{j' \in R} v_{j'i} \quad \forall i \in D \tag{26}$$

$$xd_i \leq \sum_{h \in H} dc_{ih} d_{ih} \quad \forall j \in D \tag{27}$$

$$xd_i \geq dc_i^m \sum_{h \in H} d_{ih} \quad \forall j \in D \tag{28}$$

$$\sum_{w \in W}\sum_{i \in G} d_{wi} = \sum_{w \in W}\sum_{j \in T} xt_{wj} + \sum_{j \in R} xr_j \tag{29}$$

$$\sum_{q \in Q}\sum_{h \in H} t_{qih} \leq 1 \quad \forall i \in T \tag{30}$$

$$\sum_{h \in H} t_{qih} = a_{qi} \quad \forall q \in Q, i' \in T \tag{31}$$

$$\sum_{h \in H} r_{ih} = 1 \quad \forall i \in R' \tag{32}$$

$$\sum_{h \in H} d_{ih} = 1 \quad \forall i \in D' \tag{33}$$

$$ptr_{ij} k_{ij} + ptd_{ij} z_{ij} + prd_{ij} v_{ij} \leq Pc_{ij} \quad \forall (i,j) \in A \tag{34}$$

$$\begin{aligned}
&xd_j \geq 0, z_{ij} \geq 0 && \forall j \in D, i \in T \\
&xt_{wi} \geq 0, k_{ij} \geq 0 && \forall i \in T, j \in R; w \in W \\
&xr_i \geq 0, v_{ij} \geq 0 && \forall i \in R, j \in D \\
&lo_{ik} \geq 0 && \forall i \in G, k \in K \\
&e_{ik} \geq 0 && \forall i \in (f \cup G), k \in K
\end{aligned} \tag{35}$$

$$x_{ijk} \in \{0,1\} \quad \forall i \in (f \cup G), \forall j \in (G \cup R \cup T), k \in K$$

$$r_{ih} \in \{0,1\} \quad \forall i \in R, h \in H$$

$$t_{qih} \in \{0,1\} \quad \forall i \in T, q \in Q, h \in H \tag{36}$$

$$d_{ih} \in \{0,1\} \quad \forall i \in D, h \in H$$

The objective functions (1), (2), and (3) are associated with the total cost, total risk, and $CO_2$ emissions, respectively, in the proposed HWLRPs. The first objective function is composed of seven parts, in which the first part indicates the transportation cost of waste collection. The next three parts express the transportation of waste residues. The last three parts show the fixed cost of opening recycling, treatment, and disposal facilities. The second objective function minimizes the total risk, including transportation and site risks. The transportation risk is related to the amount of waste residues transferred between facilities which are indicated by the first three terms of the objective function (2). The site risk can be described analogously to the transportation risk, except that it depends on the amount of wastes that are available at each facility. The third objective function, which is an environmental-friendly objective, is to minimize the $CO_2$ emissions from the whole hazardous waste management system.

Eq. (4) indicates that all the vehicles should start from the central depot. Eq. (5) guarantees that each vehicle leaves and arrives at the same node, which is one of the main constraints in the vehicle routing problem. Eq. (6) ensures that all generation nodes are visited once by a vehicle for collecting each type of waste. Eq. (7) ensures that the vehicles return to the central depot after they emptied their load. Eq. (8) ensures that all vehicles with collected waste can unload their waste in a treatment facility provided that their waste is compatible with the technology at that treatment facility. Eq. (9) shows that all vehicles with recyclable collected waste should empty their load at recycling facilities before going back to its origin. Eqs. (10) and (12) guarantee that the traveled distance by a heterogeneous fleet does not transgress the permissible amount. Eqs. (13)-(16) enforce the elimination of sub-tours in the problem under consideration. Eq. (17) calculates the amount of waste that is processed at each treatment facility. Eqs. (18) and (19) represent the minimum and maximum amount of waste and waste residues that are required for establishing a treatment facility.

Eq. (20) calculates the amount of waste that is processed at each treatment facility. The flow of waste residues between the treatment facilities and the recycling facilities is indicated in Eq. (21). Eqs. (22) and (23) indicate the minimum and maximum amount of waste and waste residues that are required for establishing a recycling facility. The flow of waste residues from the treatment and recycling facilities to the disposal facilities is shown in Eqs. (24) and (25). Eq. (26) measures the quantity of waste at the disposal facilities. Eqs. (27) and (28) show the minimum amount and maximum amount of waste and waste residues that are required for establishing a disposal facility. Eq. (29), which is a balanced equation, ensures that all the demand at different source nodes must be supplied. Eq. (30) states that establishing more than one treatment technology is not allowed at the transfer station. All the existing facilities are determined by Eqs. (31)-(33). Eq. (34) asserts that the transportation risk does not allow to exceed the maximum allowable tolerance capacity. Finally, Eqs. (35) and (36) impose binary and non-negative constraints for the decision variables.

*3.3. Model linearization*

The developed mathematical model has a non-linear term in the first objective function. Therefore, to solve the proposed model, the non-linear term should be transformed into linear equivalences. For this purpose, we use an exact linearization method proposed by (Azadeh et al. 2017). The non-linear term is a multiple of the binary variable $x_{ijk}$ and the continuous variable $lo_{ik}$ in the first objective function. One of the easiest ways to prevent this non-linearities is to define a new continues variable and three auxiliary constraints in the presented model. In this regard, we replace the non-linear term by the new continuous variable $xl_{ijk}$. The following three auxiliary constraints should be added to the original model to guarantee that this reformulation yields the same result as the original model:

$$xl_{ijk} \leq BM \, x_{ijk}; \qquad (i,j) \in N, \, k \in K \tag{37}$$

$$xl_{ijk} \leq lo_{ik}; \qquad (i,j) \in N, k \in K \tag{38}$$

$$xl_{ijk} \geq lo_{ik} - (1 - x_{ijk})BM; \quad (i,j) \in N, k \in K \tag{39}$$

where BM is a large number.

*3.4. Conversion to a single-objective model*

Several methods, including the ε-constraint, weighted metrics, goal programming, and lexicographic methods, are commonly applied to deal with multi-objective optimization problems. In this study, the augmented ε-constraint was introduced by (Mavrotas 2009) to deal with the proposed multi-objective optimization problem. In the augmented ε-constraint method, one of the objective functions of the problem is optimized, and the rest of the objective functions are moved to the constraints as follows:

$$Max \; \left( g_1(x) + eps \times (s_2/r_2 + s_3/r_3 + \ldots + s_p/r_p) \right)$$

*subject to*

$$g_k(x) - s_k = \varepsilon_k \quad k = 1,2,\ldots,p; \; x \in S; \; \varepsilon_k \in R^+ \tag{40}$$

in which $x$ is the vector of the decision variables, and $S$ is the solution space of the problem. $g_1(x)$, $g_2(x)$, ..., $g_p(x)$ denote the objective functions of the problem. $\varepsilon_1, \varepsilon_2, \ldots, \varepsilon_p$ are the right-hand side values of the objective functions. $r_1, r_2, \ldots, r_p$ denote the ranges of the respective objective functions, and $s_1, s_2, \ldots, s_p$ are the auxiliary variables of the relevant constraints. The value *eps* is in the interval $[10^{-6}, 10^{-3}]$.

In this method, determining the best values of $\varepsilon_k$ is critical. To find these values, the best and worst values of the objective functions considered should be obtained. To find the best value, we solve the problem considering the objective function for which we want to find the best value. For finding the worst value of an objective function, the problem is solved with the other objective functions, and the obtained values are stored. The worst value of the stored values is considered as the worst value of that objective function. By finding the best and worst values of the objective functions, an appropriate value of $\varepsilon_k$ can be determined. For this purpose, we change the value $\varepsilon_k$ between the best and worst obtained values and solve the problem. Then, the value of the first objective function at each level of $\varepsilon_k$ is analyzed to find the best one.

Regarding the above description, the multi-objective model can now be converted into an equivalent single-objective model as follows:

$$\begin{aligned} Min f_1(x) = & \sum_{i \in G} \sum_{j \in T \cup R \cup D} \sum_{k \in K} c_{ij} x_{ijk} lo_{ik} + \sum_{i \in T} \sum_{j \in D} c_{ij} z_{ij} + \sum_{i \in R} \sum_{j \in D} c_{ij} v_{ij} + \sum_{i \in T} \sum_{j \in R} c_{ij} k_{ij} \\ & + \sum_{q \in Q} \sum_{i \in T-T'} \sum_{h \in H} ft_{qih} t_{qih} + \sum_{i \in D-D'} \sum_{h \in H} fd_{ih} d_{ih} + \sum_{i \in R-R'} \sum_{h \in H} fr_{ih} r_{ih} \\ & + \left( \sum_{s \in S-\{1\}} w_s^1 \delta_s^1 \right) - \left( eps \times (S_2/r_2 + S_3/r_3) \right) \end{aligned} \tag{41}$$

*subject to*

$$\sum_{i \in T} \sum_{j \in R} Ptr_{ij} k_{ij} + \sum_{i \in T} \sum_{j \in D} Ptd_{ij} z_{ij} + \sum_{i \in R} \sum_{j \in D} Prd_{ij} v_{ij} + \sum_{w \in W} \sum_{(i,t) \in T} \sum_{h \in H} s_{th} xt_{wi}$$

$$+ \sum_{(r,i) \in R} \sum_{h \in H} s_{rh} xr_i + \sum_{(d,i) \in D} \sum_{h \in H} s_{dh} xd_i + \left( \sum_{s \in S - \{1\}} w_s^2 \delta_s^2 \right) + s_2 = \varepsilon_2 \quad (42)$$

$$\sum_{w \in W} \sum_{i \in R} Qr_{wi} xr_i + \sum_{w \in W} \sum_{i \in R} \sum_{q \in Q} Qt_{wqi} xt_{wi} + \sum_{w \in W} \sum_{i \in D} Qd_{wi} xd_i + \sum_{w \in W} \sum_{i \in T} \sum_{j \in R} QTr_{wij} dis_{ij} k_{ij}$$

$$+ \sum_{w \in W} \sum_{i \in T, R} \sum_{j \in D} QTd_{wij} dis_{ij} (z_{ij} + v_{ij}) + \left( \sum_{s \in S - \{1\}} w_s^3 \delta_s^3 \right) + s_3 = \varepsilon_3 \quad (43)$$

$$S_k \in R \quad (44)$$

Eqs. (4)-(36)

This single-objective model will be solved by the proposed solution method in the next section.

## 4. Case study

In order to validate the proposed model, a real-world case study is considered in Babol city, which is located in the northwest of the Mazandaran province, Iran. This city, with an area of 310 km² and a total population of 531930, is the most populated city in the Mazandaran province, Iran. This city is selected for the case study as it has many hospitals, clinics, and industrial companies that generate a considerable amount of hazardous waste.

To the best of our knowledge, there is no practical and strategical plan for collecting, recycling, and disposing hazardous waste in the city council. A considerable amount of generated waste is transferred to a forest area at the south of Babol, which is geographically located next to a river. Therefore, it is necessary to make a strategic plan for collecting waste and treat them in an organized way, which is the main purpose of our study. Based on the political divisions that are available on the Mazandaran province Statistics center (MPSC) website, Babol is divided into six main population districts including Central (I), Laleabad (II), Gatab (III), Bandpey-ye Gharbi (IV), Babol Kenar (V) and Bandpey-ye Shargi (VI). Fig. 2 shows the geographical map of Babol city, and its population districts are numbered from (I) to (VI). Based on (MPSC), there are 13 demand nodes within these six districts, and all districts are candidate locations for establishing facilities. In this paper, the amount of hazardous waste that is generated from industrial processes and health sectors are taken into account. According to Babol Waste Management Organization's official data for the year 2019, the amount of waste generated from industrial processes and health sectors was 131,980 tons, where 36% of them were considered as hazardous waste (BMWMC 2020).

*4.1. Data set*

In this study, all the input data we used is from the official reports of the Babol Waste Management Organization for the year 2019. The cost unit is a million Tomans (the currency of Iran, where 1 US Dollar is equal to 30000 Tomans), the distance unit is the kilometer, and the waste unit is the ton. The distances between each node and its facilities are measured by applying Google Maps. The amount of hazardous waste produced by industrial processes and health sectors in 2019 was in total 47,513 tons, and only 20% of the produced hazardous waste is transferred to facilities, and the remaining ones are disposed of improperly. Since a large amount of hazardous waste is not suitable for recycling when directly delivered from the demand nodes, a small percentage of them is sent to a recycling facility. According to the data obtained from the existing facilities, the percentage for each type of hazardous waste which are sent to corresponding facilities is shown in Table 2. Based on the research that was done by (Alumur and Kara 2007), the amount of waste residues that are sent to the recycling facility from a chemical treatment facility is 30%, and this amount for waste residues that are delivered to the recycling facility after an incineration treatment is 0% because they are only ashes. In addition, according to Alumur and Kara (2007), since the incineration process is aiming to reduce the volume of mass, a reduction of mass in this process is 80%, while this amount for a chemical treatment is 20%, since the chemical treatment aimed to reduce the hazardous characteristic of the waste. According to the data obtained from the existing facilities, after the recycling process, 5% of waste residues are sent to disposal facilities.

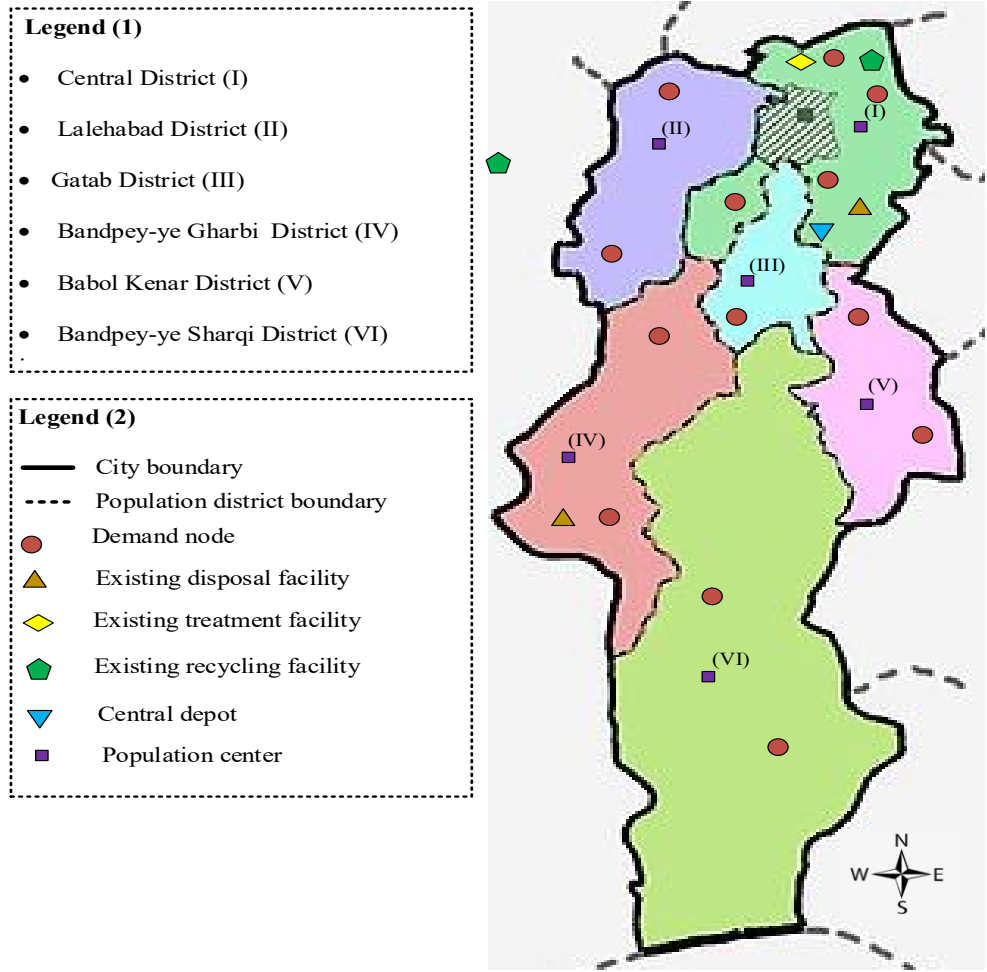

**Fig. 2.** The geographical map of Babol city.

**Table 2**
The percentage of each category of hazardous waste.

| Percentage (%) | Waste types | | | |
| --- | --- | --- | --- | --- |
| | (I) | (II) | (III) | (IV) |
| Industrial process | 15% | 13% | 20% | 5% |
| Health sectors | 10% | 30% | 15% | 8% |

The generated hazardous waste which is gathered in the demand nodes are collected by 12 trucks. Each truck collected only one type of waste because the waste should be collected separately to prevent an interaction between them (BMWMC 2020). According to the categories of the waste, the vehicles with respect to their capacities and distance limitation, deliver them to the corresponding facility. The transportation cost is calculated based on the distance that a vehicle travels to collect and deliver hazardous waste and waste residues to facilities, which are 0.01 Million Tomans per kilometer. For the establishment cost of each facility, we applied the judgment of three senior experts in the Babol Real Estate Consultants Association. The investment cost of each potential facility depends on the capacity of each facility, and it is also different for each region due to the difference in the price of the land. This is demonstrated in Table 3.

**Table 3**
The establishment cost of each potential facility (million Tomans).

| District        |    | I    | II   | III  | IV   | V    | VI   |
|-----------------|----|------|------|------|------|------|------|
|                 | 5  | 1572 | 1275 | 1463 | 1098 | 1137 | 847  |
| Capacity levels | 10 | 1932 | 1583 | 1892 | 1386 | 1408 | 1230 |
|                 | 15 | 2241 | 1987 | 2340 | 1791 | 1853 | 1596 |

According to the amount of generated waste in the selected regions, three capacity levels, including 5 (low-level),10 (medium-level),15 (high-level) tons per day, are considered for each facility. The region, where its waste production is in the range of [0-600] tons is considered to be low-level, [601-1700], and [1700-3500] are medium-level and high-level, respectively. For the region in which its waste generation is in a low-level range, the capacity level of 5 is assigned to its facilities, and the capacity level of 10 and 15 are assigned to medium-level and high-level regions, respectively. Table 4 indicates the amount of generated hazardous waste in each district. The locations of the existing facilities in each district are shown in Table 5. In addition, there are also two existing recycling facilities, including Juybar and Amol, which are located in North and West of Babol, respectively.

**Table 4**
The amount of generated hazardous waste (HW) in each region.

| District     | I      | II    | III   | IV    | V     | VI    |
|--------------|--------|-------|-------|-------|-------|-------|
| Amount of HW | 13,746 | 8,190 | 7,120 | 8,961 | 5,167 | 4,329 |

**Table 5**
The location of the existing facilities.

|                             | Types of facility |           |          |
|-----------------------------|-------------------|-----------|----------|
|                             | Recycling         | Treatment | Disposal |
| The location of the facility | (I)              | (I)       | (I), (IV)|

**Table 6**
The amount of $CO_2$ emissions from transportation and operations of the facilities.

| Types of facilities      | CO2 emission (kg emissions/t of HW) |
|--------------------------|-------------------------------------|
| Recycling                | 398                                 |
| Treatment (Incineration) | 980                                 |
| Treatment (Chemical)     | 280                                 |
| Disposal                 | 271                                 |
| Transportation           | 1.68                                |

(Liu et al. 2020)

During the operating and transportation of hazardous waste, many substances like $CO_2$ are emitted into the air. Table 6 shows the amount of produced $CO_2$ emissions by transportations and the operations of recycling, treatment, and disposal facilities according to the available data and the dependable scientific report on the footnote of Table 6. The location risk (LR) and the transportation risk (TR) of each category of hazardous waste are measured by the risk method proposed by (Nema and Gupta 1999). The values of risk consequence and risk probability are shown in Table 7 for each type of facility.

**Table 7**
The values of risk consequence and risk probability for each facility.

| Facilities | Risk | |
| --- | --- | --- |
| | Risk consequence ($\times 10^4$ people) | Risk probability ($\times 10^{-6}$) |
| Recycling | [0.01,3.32] | 20 |
| Treatment (incineration) | [0.01,3.32] | 50 |
| Treatment (chemical) | [0.01,3.32] | 60 |
| Disposal | [0.01,3.32] | 30 |

The risk consequence is defined by the number of people exposed within 800 m width of the link. The risk probability is defined as $0.4(10^{-6}/\text{km}) * 0.9 *$ link length (km), where the first term is the truck accident rate in the hazardous waste transportation on inter-city highways, the second term is the release probability given an accident of a truck with hazardous waste. The population densities for the six districts are shown in Table 8.

**Table 8**
The population density for each district.

| | District | | | | | |
| --- | --- | --- | --- | --- | --- | --- |
| | I | II | III | IV | V | VI |
| Population density (people/Km2) | 501-700 | 351-500 | 401-600 | 251-400 | 151-250 | 0-200 |

Moreover, the transportation risk is the product of three factors, including the waste risk potential, the link risk consequence, and the link risk probability. The waste risk potential is measured by an analytic hierarchy process (AHP), which is shown in Table 9 (Zhao et al. 2016).

**Table 9**
Risk potential for hazardous waste and waste residue.

| Waste type | Potential risk |
| --- | --- |
| (I) | 0.05 |
| (II) | 0.2 |
| (III) | 0.2 |
| (IV) | 0.2 |
| Disposable hazardous waste | 0.1 |
| Waste residues at a treatment facility which is recyclable | 0.05 |
| Waste residues at a treatment facility which is disposable | 0.1 |
| Waste residues at a recycling facility which is disposable | 0.1 |

*4.2. Optimal solution for the case*

The proposed single-objective model with the input parameters of the real-world case study has been solved by the CPLEX solver in the GAMS optimization software v. 24.1. The optimal solutions are illustrated in Fig. 3 and Table 10. The solution reveals that the existing facilities do not satisfy the need for the recycling, treatment, and disposal operation of the total generated waste in the studied districts. Due to the amount of generated waste, a need for a new disposal facility was clear, where according to the optimal solution, a new high-level disposal facility has been opened. As expected, a high-level facility is selected for establishing in the most populated district in which the amount of generated waste is high.

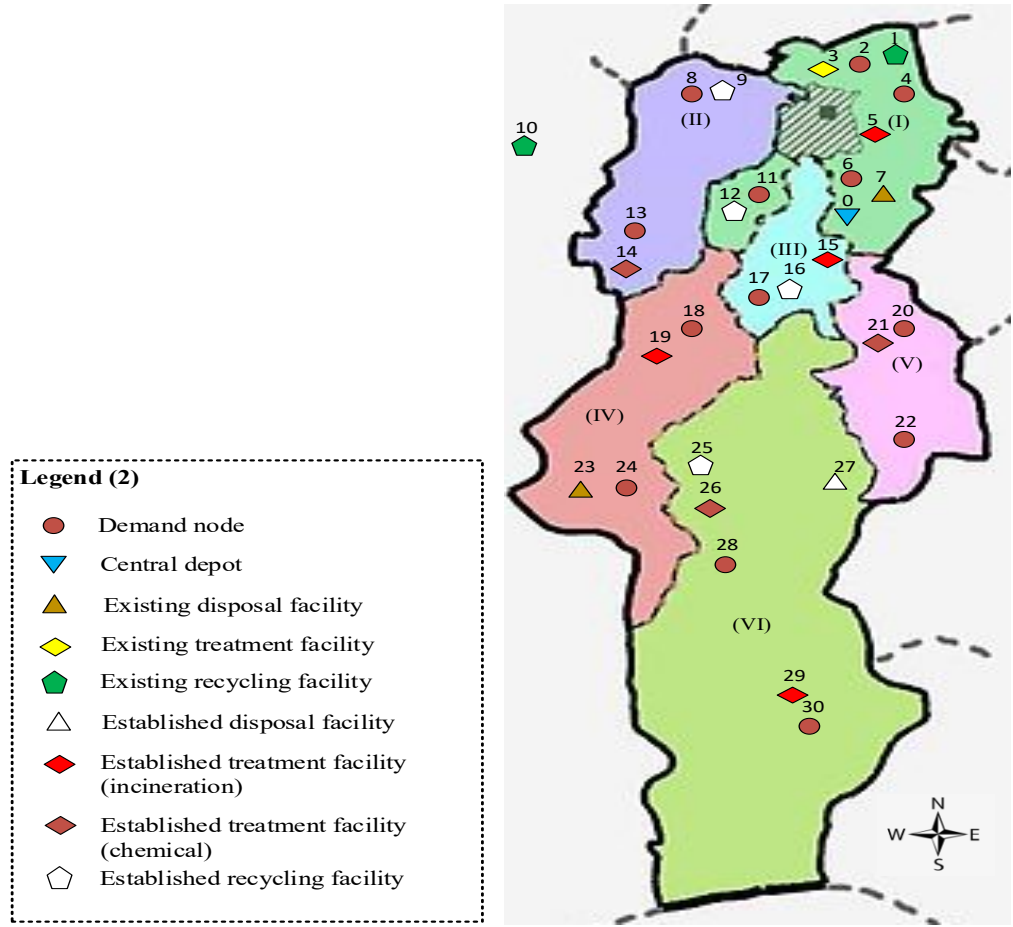

**Fig. 3.** The schematic illustration of the optimal solution for the real case study.

**Table 10**
The optimal solution for the facility location in the real case study.

| District | Facility | Number of facilities | Capacity |
|---|---|---|---|
| 1 | Recycling | 1 | High-level |
|   | Treatment (incineration) | 1 | High-level |
| 2 | Recycling | 1 | Medium-level |
|   | Treatment (chemical) | 1 | High-level |
| 3 | Recycling | 1 | Low-level |
|   | Treatment (incineration) | 1 | High-level |
| 4 | Treatment (chemical) | 1 | Medium-level |
| 5 | Treatment (chemical) | 1 | High-level |
| 6 | Recycling | 2 | Low-level, High-level |
|   | Treatment (incineration) | 1 | Low-level |
|   | Treatment (chemical) | 1 | Medium-level |
|   | Disposal | 1 | High-level |

The optimal routes of a heterogeneous fleet of vehicles for collecting the waste and transportation stage are shown in Tables 11 and 12. For representing the optimal solution of the case study in Table 11, all of the facilities in the schematic illustration are listed from 0 to 30. As one can see, Table 11 illustrates the optimal collection routes for each type of waste, where each vehicle begins from the central depot, continues its path to collect the generated waste from the demand nodes, and then it goes to the corresponding facility (recycling or treatment facility) to empty its load, and finally it returns to the central depot. For example, truck no.1 as indicated in Table 11 starts its route from the central depot, then goes to demand nodes 4, 6, 2 and 8 to collect the generated waste, and since the collected waste type is recyclable waste, it unloads its load at the recycling center 9 and terminates at the central depot. To verify the feasibility of the model, the quantity of collected waste, and length of the tour, is measured. After operating the waste at the treatment facilities, the recyclable percentage of the waste is shipped to the recycling facility, and the remaining ones are sent to the disposal facility. Moreover, at the recycling facility, the percentage of disposable waste is delivered to the disposal facility. Table 12 indicates the optimal decisions related to this stage of the problem. It is noteworthy to mention that recycling number 10, which was one of the two existing recycling facilities, was excluded from the network in the optimal solution since no amount of generated waste is delivered to this facility. In the optimal solution, due to the transportation cost, transportation and location risk, and $CO_2$ emissions, it was more sustainable to establish the new recycling facility than transporting hazardous waste to this facility. The optimal values of the objective functions are 34,481 million Tomans, 5218.4 km * people, 23381.4 tons for cost risk and $CO_2$ emissions, respectively, when separately each objective function is minimized.

**Table 11**
The optimal routes for collecting the waste.

| Waste type | Vehicle no. | Optimal route | Load's quantity (tons) | Route length (km) |
|---|---|---|---|---|
| Recyclable | 1 | 0→6→4→2→8→9→0 | 738 | 37 |
| | 2 | 0→11→13→18→17→20→16→12→0 | 952 | 48 |
| | 3 | 0→22→30→28→24→25→0 | 672 | 83 |
| Non-recyclable (incineration) | 4 | 0→6→4→2→8→5→0 | 945 | 39 |
| | 5 | 0→11→13→18→17→20→15→0 | 1063 | 42 |
| | 6 | 0→22→30→28→24→29→19→0 | 837 | 89 |
| Non-recyclable (chemical) | 7 | 0→6→4→2→8→14-0 | 884 | 49 |
| | 8 | 0→11→13→18→17→20→21→0 | 1129 | 43 |
| | 9 | 0→22→30→28→24→26-0 | 840 | 85 |
| Non-recyclable (incineration & chemical) | 10 | 0→6→4→2→8→5→0 | 652 | 39 |
| | 11 | 0→11→13→18→17→20→21→0 | 746 | 43 |
| | 12 | 0→22→30→28→24→29→19→0 | 563 | 89 |

**Table 12**
The optimal decisions for transporting the waste residues and operations at the facility centers.

| Waste residues transportation | Waste operation |
|---|---|

| Route | Load's quantity (ton) | Facility | The quantity of waste operated at each facility | Facility | The quantity of waste operated at each facility |
|---|---|---|---|---|---|
| 14→12 | 165 | 1 | 133 | 19 | 2084 |
| 21→16 | 157 | 3 | 1716 | 21 | 1632 |
| 26→25 | 99 | 5 | 2151 | 23 | 6295 |
| 14→7 | 1051 | 7 | 10341 | 25 | 1548 |
| 21→27 | 1183 | 9 | 1593 | 26 | 1842 |
| 26→23 | 1096 | 10 | - | 27 | 7297 |
| 9→7 | 721 | 12 | 1826 | 29 | 1235 |
| 12→7 | 985 | 14 | 1975 | | |
| 16→27 | 542 | 15 | 2342 | | |
| 25→23 | 939 | 16 | 1797 | | |

### 4.3. Sensitivity analysis and discussion

The results of the case study are based on the assumptions on the value of some parameters, such as the capacity level and the amount of generated waste accumulated at the demand nodes. To evaluate the sensitivity of the results with respect to these parameters, different values of these parameters were also considered as input parameters. The results and analyses are given in the following.

#### 4.3.1. The model without sustainability

To indicate the impacts of sustainability, the results of the model when sustainability is considered are compared with the results of the model without sustainability. The associated model is also solved by the CPLEX solver in the GAMS optimization software v. 24.1. In the optimal solution of the model without sustainability, three recycling centers (recycling centers 12, 16, and 25) are established.

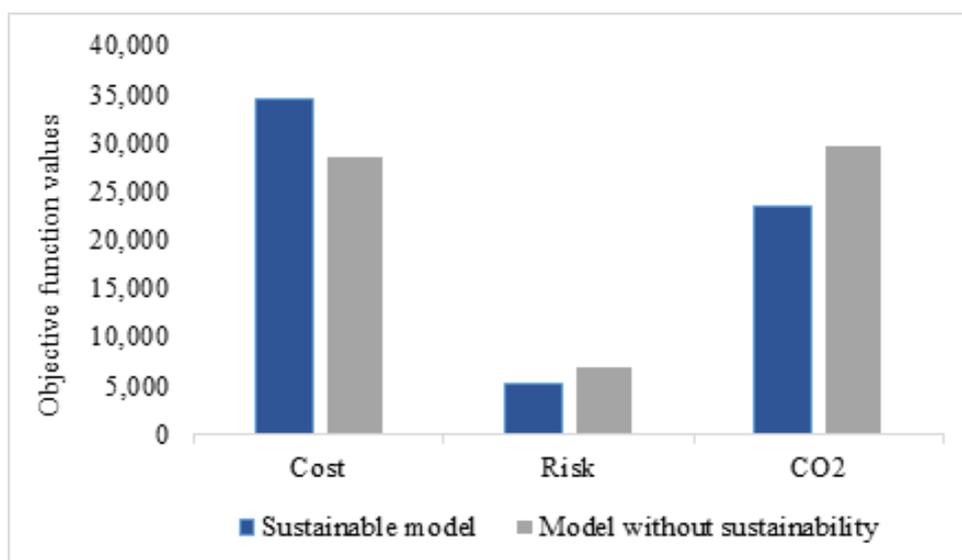

**Fig. 4.** A comparison of the objective functions.

Moreover, treatment center 15 is excluded from the optimal solution of the model without sustainability. The optimal objective function values in the model without sustainability are 28,639 million Tomans, 6831 km * people, 29571 tons, for cost risk and $CO_2$ emissions, respectively, when each objective function is individually minimized. Fig. 4 shows a comparison of the objective function of the sustainable model with the model without sustainability. The results of the sustainable model are compared to the results of the model without sustainability. The comparison indicates that, in the optimal value of the objective function, the cost of establishing facilities and transportation decreased, but the site and transportation risk along with $CO_2$ emissions increased. Since risk and $CO_2$ emissions are too critical in the hazardous waste management system in which each of them can pose a big threat to the system, it is not desirable to run a model without sustainability. These findings imply the significance of sustainability in the HWLRP.

*4.3.2. Impact of capacity level*

As mentioned earlier, according to the generated waste in each region, the facilities were established with different capacity levels. In this section, to evaluate the impact of the capacity level on the facility, three scenarios were considered associated with those decision variables which are related to the capacity level.

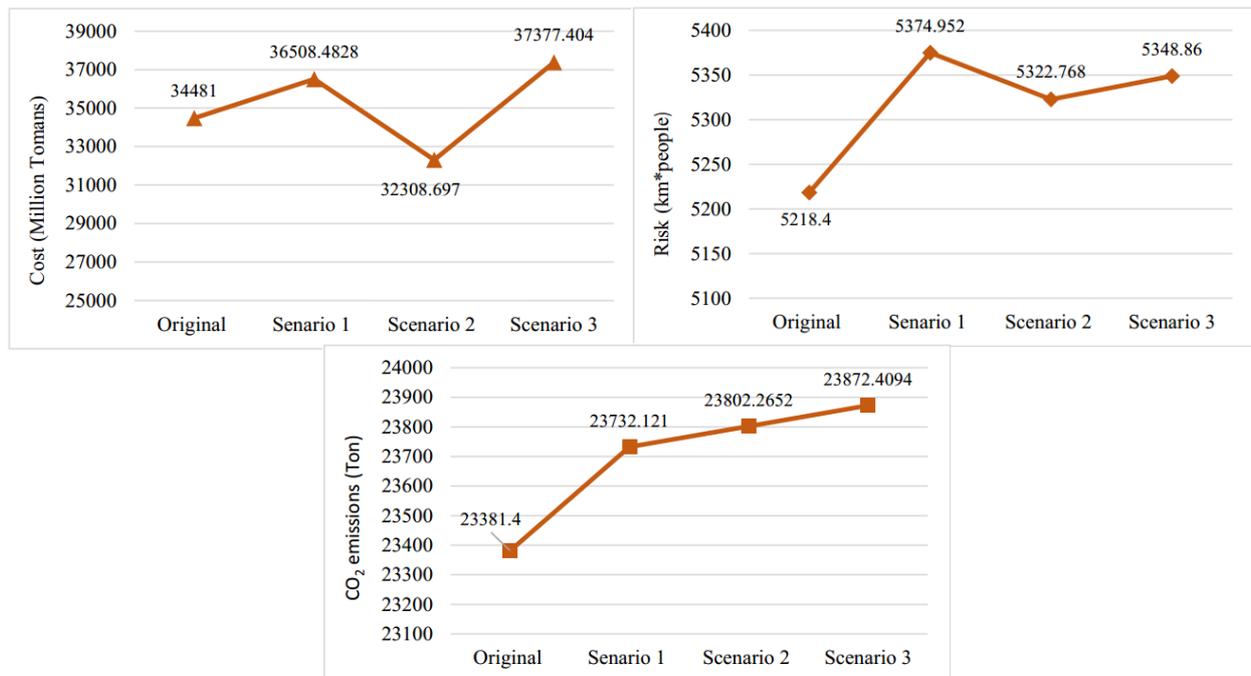

**Fig. 5.** *Sensitivity analysis of the objective functions value with respect to the changes in the capacity level*

These three scenarios are as follows: (1) remove the capacity level from the model (2) increase the capacity level (3) decrease the capacity level. Then the obtained results of solving the model associated with each scenario are compared with the results obtained in terms of the original capacity level consideration. In the first scenario, the obtained results show that under the condition that capacity planning is not considered in the problem, the total cost of the system, including investment cost, is increased by 5.88%. This is because one more facility, a recycling facility, is established in the district (V) in this case. In addition, the second (minimizing site and transportation risk) and third (minimizing $CO_2$ emissions) objective functions of the problem are increased by 3% and 1.5%, respectively. In the second scenario, by increasing the capacity level to (10, 15, 20), the total cost of the system decreased by 6.3% because recycling facility 16 and treatment facilities 15 and 26 were excluded from the optimal solution. However, the value of the risk and $CO_2$ minimization objective functions increased by 2% and 1.8%, respectively. Finally, in the last scenario in which the capacity level is decreased to (4,7,10), the obtained results indicate that under the condition of a decreasing capacity level, the total cost of the system, including investment and transportation cost, increased by 8.4%. This is mainly because by decreasing capacity level, 1 recycling and 2 treatment facilities were added to the system. Moreover, the second and third objective functions increased by 2.5% and 2.1%, respectively. Fig. 5 shows the increase and decrease of the three objective functions under three scenarios compared to the original capacity level condition. Finally, it can be concluded that considering the capacity level for the facilities would avoid imposing additional costs, and it also decreases the risk and $CO_2$ emissions, which in turn contributes to having a more sustainable system.

*4.3.3. Impact of the amount of generated waste*

In the case study, the input parameters, such as the amount of generated waste, are from the Babol Waste Management Organization's official reports.

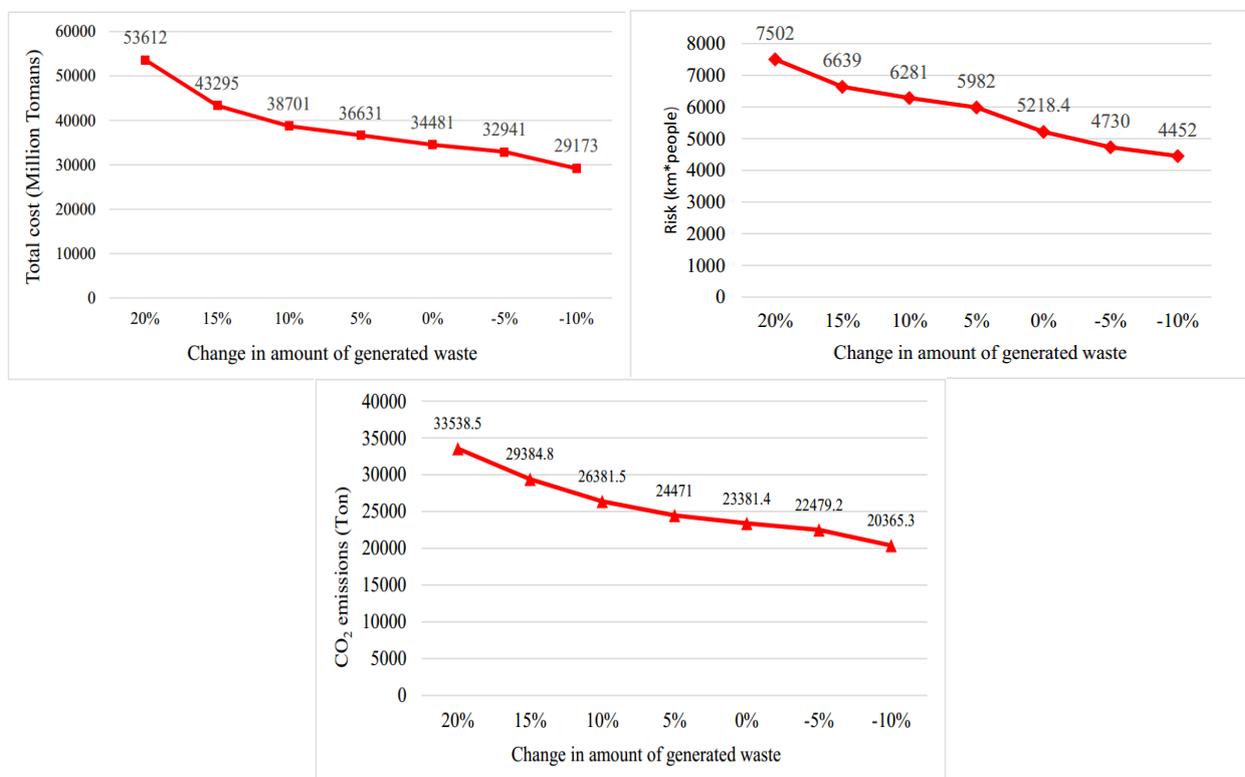

**Fig.6.** Sensitivity analysis of the objective functions value with respect to the changes in the amount of generated waste

Even though the amount of generated waste in the case study was taken from the official report of the Babol Waste Management Organization, this amount is a fixed amount. Some unprecedented happening may change the amount of generated waste in the health sector or industrial processes. For instance, these days, a widespread pandemic, COVID-19, affected the whole world, and it causes the generation of a considerable amount of infectious waste, which is characterized as hazardous waste. This shows the importance of evaluating the sensitivity of the model with respect to the amount of generated waste. In order to evaluate the sensitivity of the model with respect to changes in the amount of generated waste, the original amount was both decreased and increased by some percent to see its impact on the value of the objective function. According to the data from Fig. 6, there is a direct correlation between the change and the value of the objective functions. To be more specific, the results indicate that a minor shift in the waste amount can have a positive effect on all of the objectives.

*4.4. Managerial insights*

The purpose of this study was to design a proper HWMS with regard to the total cost, transportation, and site risk, and $CO_2$ minimization. According to the previous works in the hazardous waste management field, the lack of a strategic plan and proper financial assignment are the main concerns in Babol city. Based on the optimal results in this study, the following points would help the managers in the Babol Municipal Waste Management Center (BMWMC) to improve the HWMS.

- A lack of recycling, treatment, and disposal facilities leads to an improper disposal of the hazardous waste which poses a major threat for both environment and human. To remove these threads, many facilities at each center should be established. Fig. 3 and Table 10 present the optimal solutions which can help the managers in finding the best locations. For example, a new treatment center (treatment center 21) should be established in the district (V).

- In the existing network, most of the transportation routes are among the most populated areas to reduce the length of the routes for delivering the waste; however, it increases the associate transportation risk. It is recommended that BMWMC uses the optimal route, as shown in Table 11, for collecting the waste. For instance, truck number 1, which delivers recyclable waste, should start from the central depot, then goes to demand nodes 6, 4, 2, and 8, respectively, and unload its waste in treatment facility 5 and finally returns to the central depot.

- As an independent organization, BMWMC tries to minimize its total cost without considering the risk of $CO_2$ emissions, which causes a major problem for our environment or health. As we showed in Section 4.3.1, considering other factors like risk and emissions can lead to a better solution, even though the cost increases. Therefore, it is recommended that the government helps BMWMC with its financial aids to take these critical factors into account.

# 5. Conclusions and future work

This paper developed a mixed-integer nonlinear programming model for the sustainable hazardous waste location-routing problem. The problem is designed as a multi-objective model on a network with multiple types of hazardous waste. The proposed model determines the location of the facilities, the transportation routes for collecting the generated waste and delivering them to the facilities, and the amount of different types of waste and waste residues transported between these facilities.

To illustrate the validation of the proposed model, a real-world case study in Babol city was presented. The case study includes six districts and 13 demand nodes, where the generated waste accumulated at these nodes and then collected with a heterogeneous fleet of vehicles to transfer them to the associated facilities. The model was solved with the CPLEX solver in the GAMS optimization software v. 24.1. In the optimal solution, four new recycling facilities, seven new treatment facilities, and one new disposal facility with defined capacity levels were added to the existing network. To demonstrate the importance of sustainability, the results are compared before and after sustainability. The results showed that in the model without considering sustainability, total cost, transportation, and site risk along with $CO_2$ emissions increased, which is not desirable for a waste management system. These results emphasized the importance of sustainability in an HWMS. In addition, to evaluate the sensitivity of the model, a sensitivity analysis was performed on two parameters, including the capacity level and the amount of generated waste. The results showed that when the capacity level is removed or decreased, the values of all three objective functions increase. When the capacity level increased, the total cost decreased because fewer facilities were established. However, risk and $CO_2$ emissions increased. Moreover, there was a positive correlation between the amount of generated waste and the values of the three objective functions, which means when the amount of generated waste increases (decreases), the values of the three objective functions increase (decrease). Finally, managerial insights for hazardous waste management authorities are extracted from the final results.

For future research, this study can be extended in several directions. First, since many parameters, such as the amount of generated waste, are unknown, a stochastic version of the study is a good venue to address uncertainty. Second, future works can be conducted by developing and implementing a meta-heuristic algorithm for solving large-sized problems. Third, incorporating a time window limitation on the vehicle routes can be other features for extending the model. Finally, due to the global outbreak of the COVID-19 pandemic, which increases the amount of hazardous waste in health sectors and hospitals, considering this pandemic to design an efficient hazardous waste collection system will lessen the spread of COVID-19.

# References

Alumur, S., & Kara, B. Y. (2007). A new model for the hazardous waste location-routing problem. *Computers & Operations Research*, 34(5), 1406-1423.

Asase, M., Yanful, E. K., Mensah, M., Stanford, J., & Amponsah, S. (2009). Comparison of municipal solid waste management systems in Canada and Ghana: A case study of the cities of London, Ontario, and Kumasi, Ghana. *Waste Management*, 29(10), 2779-2786.

Asgari, N., Rajabi, M., Jamshidi, M., Khatami, M., & Farahani, R. Z. (2017). A memetic algorithm for a multi-objective obnoxious waste location-routing problem: a case study. *Annals of Operations Research*, 250(2), 279-308.

Aydemir-Karadag, A. (2018). A profit-oriented mathematical model for hazardous waste locating-routing problem. *Journal of Cleaner Production*, 202, 213-225.

Azadeh, A., Shafiee, F., Yazdanparast, R., Heydari, J., & Keshvarparast, A. (2017). Optimum integrated design of crude oil supply chain by a unique mixed integer nonlinear programming model. *Industrial & Engineering Chemistry Research*, 56(19), 5734-5746.

BMWMC, 2018. Babol manucipal waste management center. http://www.babolcity.ir.

Bronfman, A., Marianov, V., Paredes-Belmar, G., & Lüer-Villagra, A.(2016). The maxisum and maximin-maxisum HAZMAT routing problems. *Transportation Research Part E: Logistics and Transportation Review*, 93, 316-333.

Chang, N.-B., Qi, C., Islam, K., & Hossain, F. (2012). Comparisons between global warming potential and cost–benefit criteria for optimal planning of a municipal solid waste management system. *Journal of Cleaner Production*, 20(1), 1-13.

Current, J., & Ratick, S. (1995). A model to assess risk, equity and efficiency in facility location and transportation of hazardous materials. *Location Science*, 3(3), 187-201.

Darmian, S. M., Moazzeni, S., & Hvattum, L. M. (2020). Multi-objective sustainable location-districting for the collection of municipal solid waste: Two case studies. *Computers & Industrial Engineering*, 150, 106965.

Das, S., & Bhattacharyya, B. K. (2015). Optimization of municipal solid waste collection and transportation routes. *Waste Management*, 43, 9-18.

De Buck, A., Hendrix, E. M., & Schoorlemmer, H. (1999). Analysing production and environmental risks in arable farming systems: A mathematical approach. *European Journal of Operational Research*, 119(2), 416-426.

Erkut, E., & Neuman, S. (1989). Analytical models for locating undesirable facilities. *European Journal of Operational Research*, 40(3), 275-291.

Erkut, E., & Verter, V. (1998). Modeling of transport risk for hazardous materials. *Operations Research*, 46(5), 625-642.

Farrokhi-Asl, H., Tavakkoli-Moghaddam, R., Asgarian, B., & Sangari, E. (2017). Metaheuristics for a bi-objective location-routing-problem in waste collection management. *Journal of Industrial and Production Engineering*, 34(4), 239-252.

Garrido, R. A., & Bronfman, A. C. (2017). Equity and social acceptability in multiple hazardous materials routing through urban areas. *Transportation Research Part A: Policy and Practice*, 102, 244-260.

Ghaderi, A., & Burdett, R. L. (2019). An integrated location and routing approach for transporting hazardous materials in a bi-modal transportation network. *Transportation Research Part E: Logistics and Transportation Review*, *127*, 49-65.

Ghezavati, V., & Morakabatchian, S. (2015). Application of a fuzzy service level constraint for solving a multi-objective location-routing problem for the industrial hazardous wastes. *Journal of Intelligent & Fuzzy Systems*, 28(5), 2003-2013.

Goldman, A., & Dearing, P. (1975). Concepts of optimal location for partially noxious facilities. *Bulletin of the Operational Research Society of America*, 23(1), B85.

Harwood, D. W., Viner, J. G., & Russell, E. R. (1993). Procedure for developing truck accident and release rates for hazmat routing. *Journal of Transportation Engineering,* 119(2), 189-199.

Huang, G., Baetz, B., Patry, G., & Terluk, V. (1997). Capacity planning for an integrated waste management system under uncertainty: a North American case study. *Waste Management & Research*, 15(5), 523-546.

Jacobs, T. L., & Warmerdam, J. M. (1994). Simultaneous routing and siting for hazardous-waste operations. *Journal of Urban Planning and Development*, 120(3), 115-131.

Kanat, G. (2010). Municipal solid-waste management in Istanbul. *Waste Management*, 30(8-9), 1737-1745.

List, G., & Mirchandani, P. (1991). An integrated network/planar multi-objective model for routing and siting for hazardous materials and wastes. *Transportation Science*, 25(2), 146-156.

Liu, J., Huang, Z., & Wang, X. (2020). Economic and environmental assessment of carbon emissions from demolition waste based on LCA and LCC. *Sustainability*, *12*(16), 6683.

Louati, A. (2016). Modeling municipal solid waste collection: A generalized vehicle routing model with multiple transfer stations, gather sites and inhomogeneous vehicles in time windows. *Waste Management*, 52, 34-49.

Mavrotas, G. (2009). Effective implementation of the ε-constraint method in multi-objective mathematical programming problems. *Applied Mathematics and Computation*, 213(2), 455-465.

McDougall, F. R., White, P. R., Franke, M., & Hindle, P. (2008). Integrated solid waste management: a life cycle inventory: *John Wiley & Sons*.

Melachrinoudis, E., Min, H., & Wu, X. (1995). A multi-objective model for the dynamic location of landfills. *Location Science*, 3(3), 143-166.

Minh, T. T., Van Hoai, T., & Nguyet, T. T. N. (2013). A memetic algorithm for waste collection vehicle routing problem with time windows and conflicts. *Paper presented at the International Conference on Computational Science and its Applications*.

Minoglou, M., & Komilis, D. (2013). Optimizing the treatment and disposal of municipal solid wastes using mathematical programming—A case study in a Greek region. Resources, *Conservation and Recycling*, 80, 46-57.

Mohri, S. S., Asgari, N., Farahani, R. Z., Bourlakis, M., & Laker, B. (2020). Fairness in hazmat routing-scheduling: a bi-objective Stackelberg game. *Transportation Research Part E: Logistics and Transportation Review*, *140*, 102006.

MPSC, 2018. Mazandaran province Statistics center. http://www.ostan-mz.ir/.

Nema, A., & Gupta, S. (2003). Multi-objective risk analysis and optimization of regional hazardous waste management system. *Practice Periodical of Hazardous, Toxic, and Radioactive Waste Management*, 7(2), 69-77.

Nema, A. K., & Gupta, S. (1999). Optimization of regional hazardous waste management systems: an improved formulation. *Waste Management*, 19(7-8), 441-451.

Qiu, Z., Prato, T., & McCamley, F. (2001). Evaluating Environmental Risks Using Safety-First Constraints. *American Journal of Agricultural Economics*, 83(2), 402-413.

Rabbani, M., Heidari, R., Farrokhi-Asl, H., & Rahimi, N. (2018). Using metaheuristic algorithms to solve a multi-objective industrial hazardous waste location-routing problem considering incompatible waste types. *Journal of Cleaner Production*, 170, 227-241.

Rabbani, M., Heidari, R., & Yazdanparast, R. (2019). A stochastic multi-period industrial hazardous waste location-routing problem: Integrating NSGA-II and Monte Carlo simulation. *European Journal of Operational Research*, 272(3), 945-961.

Rabbani, M., Sadati, S. A., & Farrokhi-Asl, H. (2020). Incorporating location routing model and decision making techniques in industrial waste management: Application in the automotive industry. *Computers & Industrial Engineering*, 148, 106692.

Ratick, S. J., & White, A. L. (1988). A risk-sharing model for locating noxious facilities. *Environment and Planning B: Planning and Design*, 15(2), 165-179.

ReVelle, C., Cohon, J., & Shobrys, D. (1991). Simultaneous siting and routing in the disposal of hazardous wastes. *Transportation Science*, 25(2), 138-145.

Samanlioglu, F. (2013). A multi-objective mathematical model for the industrial hazardous waste location-routing problem. *European Journal of Operational Research*, 226(2), 332-340.

Stowers, C. L., & Palekar, U. S. (1993). Location models with routing considerations for a single obnoxious facility. *Transportation Science*, 27(4), 350-362.

Talyan, V., Dahiya, R., & Sreekrishnan, T. (2008). State of municipal solid waste management in Delhi, the capital of India. *Waste Management*, 28(7), 1276-1287.

Tuzkaya, G., Önüt, S., Tuzkaya, U. R., & Gülsün, B. (2008). An analytic network process approach for locating undesirable facilities: an example from Istanbul, Turkey. *Journal of Environmental management*, 88(4), 970-983.

White, L., & Lee, G. J. (2009). Operational research and sustainable development: Tackling the social dimension. *European Journal of Operational Research*, 193(3), 683-692.

Wyman, M. M., & Kuby, M. (1995). Proactive optimization of toxic waste transportation, location and technology. *Location Science*, 3(3), 167-185.

Yilmaz, O., Kara, B. Y., & Yetis, U. (2017). Hazardous waste management system design under population and environmental impact considerations. *Journal of Environmental management*, 203, 720-731.

Yu, H., Sun, X., Solvang, W. D., Laporte, G., & Lee, C. K. M. (2020a). A stochastic network design problem for hazardous waste management. *Journal of Cleaner Production*, 277, 123566.

Yu, X., Zhou, Y., & Liu, X.-F. (2020b). The two-echelon multi-objective location routing problem inspired by realistic waste collection applications: The composable model and a metaheuristic algorithm. *Applied Soft Computing*, 94, 106477.

Zahiri, B., Suresh, N. C., & de Jong, J. (2020). Resilient hazardous-materials network design under uncertainty and perishability. *Computers & Industrial Engineering*, 143, 106401.

Zhang, Y., & Zhao, J. (2011). Modeling and solution of the hazardous waste location-routing problem under uncertain conditions. In ICTE 2011 (pp. 2922-2927).

Zhao, J., Huang, L., Lee, D.-H., & Peng, Q. (2016). Improved approaches to the network design problem in regional hazardous waste management systems. *Transportation Research Part E: Logistics and Transportation Review*, 88, 52-75.

Zhao, J., & Huang, L. (2019). Multi-period network design problem in regional hazardous waste management systems. *International journal of Environmental Research and Public Health*, 16(11), 2042.

Zhao, J., & Verter, V. (2015). A bi-objective model for the used oil location-routing problem. *Computers & Operations Research*, 62, 157-168.

Zografos, K., & Samara, S. (1989). A combined location-routing model for hazardous waste transportation and disposal. *Transportation Research Record*, 1245, 52-59.